# Traveling Wave Solutions in a Reaction-Diffusion Model for Criminal Activity

H. Berestycki[a], N. Rodríguez[*][b], and L. Ryzhik[b]

[a]EHESS, CAMS, 54 Boulevard Raspail, F-75006 Paris, France.
[b]Stanford University, Department of Mathematics, Building 380, Sloan Hall Stanford, California 94305.

## Abstract

We study a reaction-diffusion system of partial differential equations, which can be taken to be a basic model for criminal activity, first introduced in [3]. We show that the assumption of a populations natural tendency towards crime significantly changes the long-time behavior of criminal activity patterns. Under the right assumptions on these natural tendencies we first show that there exists traveling wave solutions connecting zones with no criminal activity and zones with high criminal activity, known as *hotspots*. This corresponds to an invasion of criminal activity onto all space. Second, we study the problem of preventing such invasions by employing a finite number of resources that reduce the payoff committing a crime in a finite region. We make the concept of wave propagation mathematically rigorous in this situation by proving the existence of entire solutions that approach traveling waves as time approaches negative infinity. Furthermore, we characterize the minimum amount of resources necessary to prevent the invasion in the case when prevention is possible. Finally, we apply our theory to what is commonly known as the *gap* problem in the excitable media literature, proving existing conjectures in the literature.

## 1 Introduction

The work of Short and collaborators on the mathematical modeling of criminal behavior in [19] has generated much interest, not only in the mathematical community, but in other academic and non-academic spheres. In [19] the authors introduced a parabolic system of PDEs to study the dynamics of crime *hotspots*, which are spatial and temporal areas of high density of crime. This work has been followed by a vast amount of research geared towards using mathematical tools in an effort to gain some understanding of the phenomenon of crime, see for example [14, 15]. This includes the work of Berestycki and Nadal in [3], where they introduce a more general class of reaction-diffusion system of equations. In this work the authors introduce the concept of a *warm-spot*, which are areas with some criminal activity, but not enough to be classified as a hotspot. The use of PDE systems to model this phenomenon provides a different, interesting, and insightful perspective. Furthermore, a PDE model allows for additional mathematical tools not available in more commonly used discrete models. The focus of this work is the study a class of reaction diffusion systems of partial differential equations introduced in [3], as means of continuing to explore criminal behavior. In particular, we are interested in the role that population's *natural tendencies towards criminal activity* affect the long time patterns. Our

---

[*]corresponding author: nrodriguez@math.stanford.edu

starting point is the system introduced in [3]:

$$s_t(x,t) = \Delta s(x,t) - s(x,t) + s_b(x) + (\rho(x) - c(x,t))u(x,t) \tag{1a}$$

$$u_t(x,t) = \frac{1}{\tau_u}\left(\Lambda(s) - u(x,t)\right) \tag{1b}$$

$$c_t(x,t) = \frac{1}{\tau_c}\left(-c(x,t) + c_0 + \eta_c c_0 \left(\frac{u(x,t)p(x)}{\int u(x,t)p(x)} - 1\right)\right), \tag{1c}$$

defined for $x \in \mathbb{R}^n$ and $t \geq 0$. The unknowns are: the *moving average of crime*, $u(x,t)$, the population's *propensity to commit a crime*, $s(x,t)$, and the *cost of committing a crime*, $c(x,t)$. This model, as the model studied in [19], is based on two fundamental assumptions from criminology theory: *routine activity theory* [6] and the *repeat and near-repeat victimization effect*. Routine activity theory states that the two essential elements for a crime to occur are a motivated agent and an opportunity, and neglects other factors. This has the additional benefit of making the problem mathematically tractable. The repeat and near-repeat victimization effect, which has been observed in real crime data, is the effect that crime in an area leads to more crime [12, 21].

We now briefly explain how (1) incorporates the above assumptions, see [3] for more details. Let us start with (1b). The most fundamental assumption is that the moving average of crime $u(x,t)$ and the propensity to commit crime $s(x,t)$ positively influence each other. If $s(x,t)$ is non-positive then it does not affect the crime average but if $s(x,t)$ becomes positive, its influence on the criminal activity increases. The fraction of the population at a given place and time that will commit a crime is given by

$$\Lambda(s) = \begin{cases} 0 & s < 0 \\ 1 - e^{-\beta s} & s \geq 0. \end{cases} \tag{2}$$

Here $\beta > 0$ is a given constant. When $s(x,t) = 0$ the crime will disappear on the time scale $\tau_u$, and the balance of these two effects gives (1b).

The propensity to commit a crime, $s(x,t)$, evolves proportionally to the amount of crime: it is observed in the crime data that crime leads to more crime, this is the 'repeat and near-repeat victimization effect'. Hence, we assume that $s(x,t)$ increases proportionally to $u(x,t)$ with a rate that is the difference of two effects: $\rho(x)$ measures the payoff of committing a crime and $c(x,t)$ is its cost. Note that $\rho(x)$ is naturally spatially heterogeneous as different neighborhoods have different payoff for a successful crime. A good example to keep in mind is residential burglaries, where the net payoff of a successful burglary is higher in wealthier neighborhoods. In addition, the diffusion term in (1a) takes into account the non-local effect that can also be seen on the example of residential burglaries. The function $s_b(x)$ in (1a) is the equilibrium value, or a measure of the population innate tendency towards criminal activity. If $s_b < 0$ the population is assumed to have a *natural anti-crime tendency*, $s_b = 0$ assumes a *natural indifference towards criminal activity*, and $s_b > 0$ assumes an *natural tendency towards criminal activity*. Furthermore, it corresponds to the base attractiveness value $A_o(x)$ in the model introduced in [19]. To the authors' knowledge there have been no studies of the important effect that the base attractiveness value has on the long time behavior of the solutions. Here we show that $s_b(x)$ plays a significant role in the long term behavior of criminal activity patterns.

Equation (1c) models what is referred to in [3] as an *adaptive cost*, assuming that resources are allocated at high crime rate areas. This is known as *hotspot policing*. As a first step, we study a model that does not consider the dynamics of the cost (1)

$$s_t(x,t) = \Delta s(x,t) - s(x,t) + s_b + \alpha(x)u(x,t) \tag{3a}$$

$$u_t(x,t) = \Lambda(s) - u(x,t). \tag{3b}$$

Furthermore, we assume that $s_b$ is spatially homogeneous ($s_b(x) \equiv s_b$). We plan to address the full model in the near future.

We consider initial conditions that only depend on one spatial variable, $x_1$,

$$s(x,0) = s_0(x_1) \tag{4a}$$
$$u(x,0) = u_0(x_1). \tag{4b}$$

When $\alpha < 0$ the system decays to equilibrium, and the dynamics is not particularly interesting, so we assume that $\alpha(x) \geq 0$, which is also the more realistic case sociologically.

A steady state solution satisfies

$$u = \Lambda(s) \quad \text{and} \quad s_{xx} + f(x,s) = 0 \tag{5}$$

where

$$f(x,s) = -s + s_b + \alpha(x)\Lambda(s). \tag{6}$$

We see that indeed the sign of $s_b$ determines the number of solutions to (5) and thus affects the long time behavior of solutions. Equation (5) provides a natural connection between the parabolic system (3) and the single parabolic equation when the reaction term is homogeneous ($\alpha$ is constant)

$$s_t = s_{xx} + f(s). \tag{7}$$

The parabolic equation (7) has been the subject of much study since its introduction by Fisher in [9] as a model for the spread of advantageous genetic traits in a population. Since then (7) has been used to model numerous other phenomena, from flame propagation [1] to nerve pulse propagation [16].

We consider three aspects of (3): first, the invasion of high criminal activity into areas with low initial criminal activity – this is related to the existence of traveling wave solutions. Second, we study the prevention of such invasions through what is commonly known, in the excitable media literature, as the *gap problem* [13]. In this formulation, limited resources are employed in a finite region in order to reduce the strength of reaction term (which fosters propagation). We prove that, in the right parameter regime, it is possible to prevent the propagation of criminal activity. Furthermore, we characterize the minimal length of the interval, where the reaction term is replaced by decay, required to prevent the crime wave propagation. We also prove the existence of a unique entire in time solution that approaches a traveling wave solution in the limit $t \to -\infty$. This result rigorously defines the notion of wave propagation in the heterogenous problem. Finally, we address the issue of splitting resources and prove that splitting resources is never beneficial.

Previous work on preventing propagation of traveling waves in excitable media includes, for instance, applications in physiology see [17] and in chemical reaction theory [5]. Of particular relevance is the work of Lewis and Keener [13], who consider the gap problem for the single parabolic bistable equation (7) in $\mathbb{R}$. Through a bifurcation analysis and a geometric approach they prove the existence of a unique critical length that stops the propagation of the wave. We generalize the results of [13] to a system of equations using different methods that seem to us to be more adaptable to other problems. In addition, we prove that there is a unique entire in time solution of the gap problem that converges to a traveling wave as $t \to -\infty$. This result also applies to (7).

*Outline:* We describe the main results in Section 2. The bistable system is treated in Section 3 and the monostable system is treated in Section 4. We devote Section 5 to the proof of existence of an entire solution to (7), which approaches a traveling wave solution as time approaches negative infinity. We conclude with a discussion in Section 6.

**Acknowledgment.** N. Rodríguez would like to thank Jacob Bedrossian for his helpful discussion. N. Rodríguez was partially supported by the NSF Postdoctoral Fellowship in Mathematical Sciences DMS-1103765. The research leading to these results has received funding from the European Research Council under the European Union's Seventh Framework Programme (FP/2007-2013) / ERC Grant Agreement n.321186 - ReaDi -Reaction-Diffusion Equations, Propagation and Modelling. Part of this work was done while Henri Berestycki was visiting the University of Chicago. He was supported by an NSF FRG grant DMS-1065979. L. Ryzhik was supported by the NSF grant DMS-0908507.

## 2 Main results

### 2.1 Definitions, Notation, and Preliminaries

In order to deal with non-negativity of solutions, we rewrite (3) in terms of $\tilde{s} = s - s_b$:

$$\tilde{s}_t(x,t) = \partial_{xx}\tilde{s}(x,t) - \tilde{s}(x,t) + \alpha_L(x)u(x,t) \tag{8a}$$
$$u_t(x,t) = g(\tilde{s}) - u(x,t), \tag{8b}$$

where, $s_b \in \mathbb{R}$, $g(\tilde{s}) = \Lambda(\tilde{s} - s_b)$ and

$$\alpha_L(x) = \begin{cases} \alpha & x < 0 \text{ and } x > L \\ 0 & 0 < x < L, \end{cases} \tag{9}$$

for $L \geq 0$ and $\alpha > 0$. The case $L = 0$ corresponds to the problem on a *homogeneous medium* and the case when $L > 0$ corresponds to that in a *heterogeneous medium*. The steady state version of (8) is

$$u = g(s) \tag{10a}$$
$$s_{xx} = -f_L(x,s), \tag{10b}$$

with

$$f_L(x,s) = -s + \alpha_L(x)g(s). \tag{11}$$

Here we have replaced $\tilde{s}$ with $s$ for notational simplicity. Note that $L > 0$ corresponds to the gap problem (7) where the reaction term is replaced by the decay term $-s$ in an interval of length $L$.

When $L = 0$ for simplicity we denote the reaction term by $f(s)$, which has either one, two, or three zeros depending on the parameters $s_b, \alpha, \beta$. In the *bistable case* $f$ has three zeros, and we assume without loss of generality that $s_b$, $\beta$, and $\alpha$ are such that there exists a unique $a \in (0,1)$ so that

$$f(0) = f(a) = f(1) = 0, \quad 0 < s_b < a < 1, \quad \text{and} \quad f'(0) < 0, \ f'(a) > 0, \ f'(1) < 0. \tag{12}$$

In this case, the steady state $(0,0)$ represents zero criminal activity, $(a, g(a))$ is a warm-spot, and $(1, g(1))$ is a hotspot. The states $(0,0)$ and $(1, g(1))$ are stable, and $(a, g(a))$ is unstable. Thus, in the long term we expect a situation where there is either a large amount of criminal activity or a complete lack of it.

On the other hand, the system is *monostable* if the parameters are chosen so that

$$f(0) = f(1) = 0, \quad f'(0) > 0, \ f(u) > 0, \text{ for } u \in (0,1). \tag{13}$$

For the monostable case we have that $(0,0)$ is unstable, as $\lim_{s \to 0^+} \alpha g'(s) > 1$ (which is required in order for $f(s)$ to be monostable), and $(1, g(1))$ is stable. The case when there is a single steady state solution is the least interesting and we do not consider it here. Finally, for technical reasons we assume that $\beta$ is large enough such that

$$\frac{\log \alpha \beta}{\beta} < 1. \tag{14}$$

We first prove that there exists a traveling wave solution of (8) that connects two steady states: a hotspot invading zero-criminal activity region. Let us recall [8] that for the single parabolic equation (7) there exists a unique monotone traveling wave front solution $S(x - \tilde{c}t)$ to (7), which satisfies

$$S''(z) + \tilde{c}S'(z) + f(S) = 0, \quad \lim_{x \to -\infty} S(z) = 1, \quad \lim_{x \to \infty} S(z) = 0.$$

Furthermore, $S'(z) < 0$ for $|z| < \infty$ and $\tilde{c} > 0$ iff

$$\int_0^1 f(s)\, ds > 0. \tag{15}$$

The monostable case differs as traveling wave solutions exists for a semi-infinite interval of speeds. In particular, there exists a $c^\star$ such that there are traveling wave solutions $U(x - ct)$ for any $c \geq c^\star$. Here, it is always the case that steady state $U = 1$ invades $U = 0$.

We have the following comparison principle.

**Lemma 1** (Comparison Principle). *Let $L \geq 0$, $(s_1, u_1)$ and $(s_2, u_2)$ be two solutions of (8) such that $s_1(x, 0) \geq s_2(x, 0)$ and $u_1(x, 0) \geq u_2(x, 0)$. Then $s_1(x, t) \geq s_2(x, t)$, $u_1(x, t) \geq u_2(x, t)$ for all $t > 0$.*

*Proof.* Let $s = s_1 - s_2$ and $u = u_1 - u_2$, then $u$ and $s$ satisfy

$$s_t = s'' - s + \alpha_L(x)u$$
$$u_t = \frac{g(s_1) - g(s_2)}{s_1 - s_2} s - u.$$

Since $g(s)$ is monotone increasing we have that

$$\frac{g(s_1) - g(s_2)}{s_1 - s_2} \geq 0.$$

Moreover, $\alpha_L(x)$ is also positive so by the standard monotone iteration scheme for the construction of solutions one obtains that $s(x, t), u(x, t) \geq 0$ for all $t > 0$ and we conclude. □

**Definition 1** (Supersolution). A function $\psi_+(x) \in C^2(\mathbb{R}\setminus \{0, L\})$ is a *supersolution* to (10) if it satisfies $\partial_{xx}\psi_+ + f(\psi_+) \leq 0$ and for $\xi \in \{0, L\}$

$$\lim_{x \to \xi^+} \psi'_+(x) \leq \lim_{x \to \xi^-} \psi'_+(x). \tag{16}$$

A *subsolution*, $\psi_-(x)$, is defined similarly by reversing the inequalities above. The above definition is easily extended to a finite number of discontinuities.

## 2.2 Summary of Results

In this subsection we summarize our results. Our main interest is the case when a population is said to have *anti-criminal tendencies*, which mathematically corresponds to the case the system (8) is bistable and we discuss this case first.

### 2.2.1 Bistable System

We first study the case of a homogeneous medium (*i.e* $L = 0$). We show that there exist traveling wave solutions, $(\psi, \phi, c)$ with $c \in \mathbb{R}^+$, such that for all $z \in \mathbb{R}$

$$\begin{cases} \psi''(z) + c\psi'(z) - \psi(z) + \alpha\phi(z) = 0 \\ c\phi'(z) + g(\psi(z)) - \phi(z) = 0 \\ 0 \leq \psi(z), \phi(z) \leq 1 \\ \psi(-\infty) = 1,\ \phi(-\infty) = g(1),\ \psi(+\infty) = 0,\ \phi(+\infty) = 0. \end{cases} \tag{17}$$

**Theorem 1** (Traveling Wave Solution for the Bistable System). *Let $L = 0$ and $f(s) = -s + \alpha g(s)$ be bistable, that is $f(s)$ satisfies (12). Then there exists a unique $(\psi(z), \phi(z), c^*)$ with $c^* \in \mathbb{R}$, up to translations, to (17), Furthermore,*

(i) *$c^* > 0$ iff (15) is satisfied.*

(ii) *$\psi'(z) < 0$ and $\phi'(z) < 0$ for all $|z| < \infty$.*

Now, consider initial data that roughly has the shape of a traveling wave. In particular,

$$\liminf_{z \to -\infty} s_0(z) > a, \quad \limsup_{z \to \infty} s_0(z) < a \tag{18a}$$

$$\liminf_{z \to -\infty} u_0(z) > g(a), \quad \limsup_{z \to \infty} u_0(z) < g(a). \tag{18b}$$

The next result states that initial conditions that satisfy (18) lead to the eventual invasion of criminal activity everywhere.

**Theorem 2** (Homogeneous Case: Invasion of Criminal Activity and Stability). *Let $L = 0$ and $(s(x,t), u(x,t))$ be solutions to (8) with initial conditions, $(s_0(x), u_0(x))$, such that $0 \leq u_0(x) \leq g(1)$ and $0 \leq s_0(x) \leq 1$ satisfying (18). Furthermore, let $(\psi(z), \phi(z), c)$ be the traveling wave solutions. Then there exists $z_1, z_2 \in \mathbb{R}$ and bounded, positive functions $q_u(t)$ and $q_s(t)$ that approach zero as $t \to \infty$ such that*

$$\psi(z + z_1) - q_s(t) \leq s(z + ct, t) \leq \psi(z + z_2) + q_s(t) \tag{19a}$$

$$\phi(z + z_1) - q_u(t) \leq u(z + ct, t) \leq \phi(z + z_2) + q_u(t), \tag{19b}$$

*for all $t \geq 0$, $z \in \mathbb{R}$. If additionally, the initial conditions satisfy:*

$$|u_0(z) - \phi(z + z_0)| \leq \epsilon \quad \text{and} \quad |s_0(z) - \psi(z + z_0)| \leq \epsilon \tag{20}$$

*for some $z_0, \epsilon > 0$. Then there exists a function, $\omega(\epsilon)$, which satisfies $\lim_{\epsilon \to 0} \omega(\epsilon) = 0$ and*

$$|u(z + ct, t) - \phi(z + z_0)| \leq \omega(\epsilon) \quad \text{and} \quad |s(z + ct, t) - \psi(z + z_0)| \leq \omega(\epsilon), \tag{21}$$

*for all $t \geq 0$, $z \in \mathbb{R}$. In particular,*

$$\lim_{x, t \to \infty} u(x, t) = g(1).$$

The inequality (21) tells us that once a solution is close to a traveling wave solution it will remain close for all time. Numerical simulations indicate that the solutions of the initial value problem converge exponentially fast to traveling waves. The next result is for the heterogeneous problem ((8) with $L > 0$) and a bistable reaction term.

**Theorem 3** (Unique Entire Solutions). *There exists a unique entire solution (up to a time shift) $(\bar{s}(x,t), \bar{u}(x,t))$ to (8) with $f(x,s)$ bistable, defined on $\mathbb{R} \times \mathbb{R}$, such that $0 < \bar{s}(x,t) < 1$, $0 < \bar{u}(x,t) < g(1)$, $\bar{s}_t(x,t) > 0$, and $\bar{u}_t(x,t) > 0$ for all $(x,t) \in \mathbb{R} \times \mathbb{R}$. Furthermore, if $(\psi(z), \phi(z), c)$ are the unique solutions from Theorem 1 then*

$$\bar{s}(x,t) \to \psi(x - ct), \quad \bar{u}(x,t) \to \phi(x - ct), \tag{22}$$

*as $t \to -\infty$ uniformly in $x \in \mathbb{R}$.*

Such solutions are called transition fronts, and one interpretation for the prevention of the invasion of criminal activity is that the transition waves tend to a steady state and do not propagate as $t \to +\infty$.

The comparison principle for (8) enables us to reduce the problem of blocking the generalized traveling fronts, $(\bar{s}(x,t), \bar{u}(x,t))$, to finding a steady state solution, $s_L(x)$, that is monotonically decreasing. In particular, we require that the steady state solution satisfy

$$\lim_{x \to -\infty} s_L(x) = 1 \quad \lim_{x \to -\infty} s'_L(x) = 0. \tag{23}$$

The condition at negative infinity allows for the initial propagation of the wave and the condition at positive infinity guarantees that the transition fronts are blocked. Thus, the problem is

$$s'' + f_L(x, s) = 0 \tag{24}$$

where,

$$f_L(x, s) = \begin{cases} f(s) & x < 0 \text{ and } x > L \\ -s & 0 \leq x \leq L, \end{cases} \tag{25}$$

with the boundary condition (23).

Our next result, and the crux of this work, is on the existence of critical gap length $L^\star$, such that for $L \geq L^\star$, transition fronts are blocked. On the other hand, for $L < L^\star$ the fronts are delayed, but eventually propagate (see Figure 1a). Before stating this result, we define $b \in (0,1)$ such that $\mathcal{F}(b) = 0$, where

$$\mathcal{F}(s) := \int_0^s f(v) \, dv.$$

**Theorem 4** (Blocking the Generalized Traveling Fronts). *Let $(\bar{s}(x,t), \bar{u}(x,t))$ be solution to (8) from Theorem 3, then there exist a unique $L^\star > 0$ such that*

(i) *if $L \geq L^\star$ then the propagation of the waves is blocked, i.e there exists a solution $s_L(x)$ to (24) that satisfies (23).*

(ii) *if $L < L^\star$ then the waves propagate. That is, for any $\epsilon > 0$ there exists $x_\epsilon > L$ and $t_\epsilon > 0$, such that*

$$\bar{s}(x,t) > 1 - \epsilon \quad \text{for} \quad x > x_\epsilon, \ t \geq t_\epsilon \tag{26a}$$
$$\bar{u}(x,t) > g(1) - \epsilon \quad \text{for} \quad x > x_\epsilon, \ t \geq t_\epsilon. \tag{26b}$$

*Furthermore, $L^\star$ is characterized as the unique value such that (24), with $L = L^\star$, has a unique monotone decreasing solution, $s_{L^\star}(x)$, satisfying $s_{L^\star}(L^\star) = b$ and $s'_{L^\star}(L^\star) = 0$.*

The results of Theorem 4 are in agreement with numerical simulations, see Figure 1, which shows the evolution of two numerical solutions with the only difference being the lengths of the gap. On the left, $L = .84$ and $t = 350$, one observes that after some time the wave propagates. Furthermore, the traveling wave front is eventually reformed after some time. On the right, $L = .85$ and $t = 800$ and, contrary to the previous case, the invasion of criminal activity is prevented.

Finally, we address the question of splitting resources for the bistable system. Given that the minimum amount of resources required to stop the propagation of the generalized wave is $L = L^\star$, we now ask the following question: if we split the resources into two regions with decay, that is two intervals with lengths $L_1$ and $L_2$, which are separated by a distance $d$, does there exist a split and a distance $d$ between these two intervals such that $L_1 + L_2 < L^\star$ is sufficient to stop the propagation of the waves? See Figure 2 for an illustration of the problem. The next result states that splitting resources can never help minimize the total amount of resources used.

**Proposition 1** (Splitting Resources). *Let $L < L^\star$, with $L^\star$ from Theorem 4. Consider an arbitrary split of an interval of length $L$ into two sub-intervals of length $L_1$ and $L_2$ (as described above). Then, for any distance $d \geq 0$ between the two intervals of length $L_1$ and $L_2$ the generalized traveling waves will propagate.*

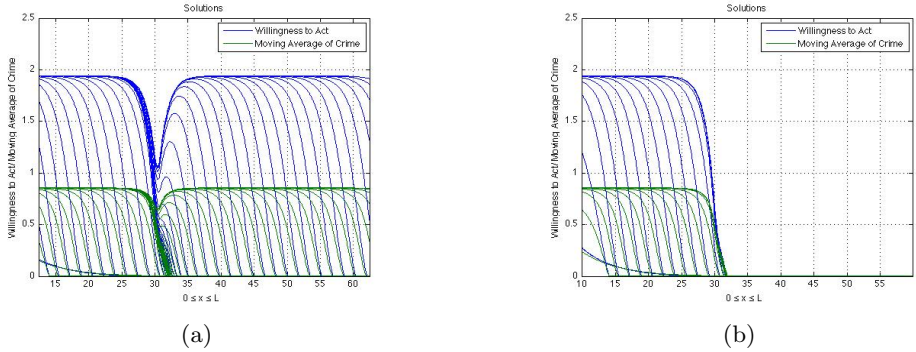

(a)            (b)

Figure 1: Two simulations with the same parameters but with different values of $L$, both where ran with a mesh-size of $\Delta x = .01$. On the right $(b)$ shows the first numerical value that blocks the wave, $L^\star = .85$, on the left $(a)$ the simulation was ran with $L = .84$.

#### 2.2.2 Monostable System

In this subsection we discuss the case when the population's natural tendencies towards crime is *indifference* (this corresponds to $s_b = 0$). Mathematically, this corresponds to a monostable system. We first state the corresponding result for traveling wave solutions in the monostable case.

**Theorem 5** (Traveling Wave Solution for the Monostable System). *Let $L = 0$ and $f(s)$ be monostable, there exists a $c^* \in \mathbb{R}^+$ and solutions $(\psi_c(z), \phi_c(z))$ with $z = x - ct$ to (17) for any $c \geq c^*$. Furthermore, when $c < c^*$, such waves do not exist.*

An interesting conclusion is that assuming that $s_b < 0$ was crucial in preventing the propagation of crime. In fact, when the population is *neutral* towards criminal activity, then no amount of limited resources can prevent the propagation of crime waves.

**Theorem 6** (Monostable Criminal Activity). *Let $L > 0$ and $(s(x,t), u(x,t))$ be solutions to (8) with $f(s)$ monostable, then the waves propagate. In particular, (26) holds for some $x_\epsilon > L$ and $t_\epsilon > 0$.*

#### 2.2.3 Nerve-Pulse Propagation Equation

As noted earlier, system (8) is related to equation (7). In fact, our methods and results carry over to the bistable (7), and the results of [13] are included in Theorem 4. We can also qualify the invasion process by proving a similar result to that in Theorem 3, adding to the picture of the gap problem for the single parabolic bistable equation.

**Theorem 7** (Unique Entire Solution). *Let $L > 0$, there exists a unique entire solution $\bar{s}(x,t)$ to (7) on $\mathbb{R} \times \mathbb{R}$ with $0 < \bar{s}(x,t) < 1$, $\bar{s}_t(x,t) > 0$ for all $(x,t) \in \mathbb{R} \times \mathbb{R}$. Furthermore,*

$$\bar{s}(x,t) - \psi(x + ct) \to 0, \tag{27}$$

*as $t \to -\infty$.*

## 3   Bistable Case: Anti-Criminal Activity Tendencies

This section is devoted to proving the results for the bistable case.

## 3.1 On Traveling Wave Solutions for the Homogeneous Problem

*Proof.* (Theorem 1) The existence and monotonicity of the traveling wave solution, $(\phi, \psi, c)$, which satisfies (17), follows from the general result in [20] (Theorem 2.1) on the existence of traveling waves for monotonic systems. It is straightforward to verify that the conditions of this theorem are satisfied. To prove part *(i)* we let $S(z)$ be the traveling wave solution of (7) with unique speed $\tilde{c}$. Recall that it satisfies

$$S'' + \tilde{c}S' + f(S) = 0. \tag{28}$$

Consider the moving coordinates $z = x - ct$ and define $w(z,t) = s(z+ct, t)$ and $v(z,t) = u(z+ct, t)$. Then any solution $\mathbf{v} = (w, v)$ to (8), satisfies $N\mathbf{v} = 0$ with

$$N\mathbf{v} = \begin{cases} w_t - w_{zz} - cw_z + w - \alpha v \\ v_t - cv_z - g(w) + v. \end{cases} \tag{29}$$

Let $\tilde{c} > 0$ and assume for contradiction that $c \leq 0$. Applying the operator $N$ to $\mathbf{v_1} = (S(z), g(S(z)))$ gives

$$N\mathbf{v_1} = \begin{cases} -S''(z) - cS'(z) - f(S(z)) \\ -cg'(S(z))S'(z) - g(S(z)) + g(S(z)). \end{cases}$$

Given that $S'(z) < 0$ and that $S(z)$ satisfies (28) we see that $-S'' - f(s) < 0$. Also, $-cS'(z) \leq 0$. Thus $N\mathbf{v_1} \leq 0$, which implies that $\mathbf{v_1}$ is a subsolution. We just need to check that there are not issues at infinity. To do this, recall that the traveling waves are unique modulo translations. We take $z^\star$ large enough so that

$$S(z + z^\star) < \psi(z) \text{ and } g(S(z + z^\star)) < \phi(z) \text{ for } z \leq R, \tag{30}$$

for some real number $R$ and

$$S(z + z^\star) < \delta, \text{ for } z > R,$$

for some small $\delta$.

We show that (30) actually holds for all $z \in \mathbb{R}$. Recall that $\tilde{c} \leq 0$ and $S'(z), \psi'(z), \phi'(z) < 0$ and that we assumed that $c < 0$. Thus,

$$\psi''(z) + f(\psi(z)) \leq 0$$
$$S''(z) + f(S(z)) \geq 0.$$

Thus, $w = S - \psi$ satisfies

$$w''(z) + q(x)w \geq 0,$$

with

$$q(x) = \frac{f(S) - f(\psi)}{S - \psi}.$$

Note that $q(x) \leq 0$ if $w \geq 0$, as $S \leq \delta$. We apply the strong maximum principle to the set $\{w > 0\}$ and conclude that $w \leq 0$ since $w = 0$ on the boundary, which is a contradiction. Thus, $w \leq 0$ and so $S(z)$ is below $\psi(z)$ for all $z \in \mathbb{R}$. A similar argument works for $g(S(z))$ and $\phi(z)$. Finally, since $S$ and $g(S)$ are traveling waves moving to the right and $c \leq 0$ implies that the traveling wave solutions $(\psi, \phi)$ are either moving to the left or stationary (if $c = 0$), this leads to a contradiction. The same argument works when $\tilde{c} \leq 0$ and this proves *(i)*. □

## 3.2 On the Dynamics of the Bistable Homogeneous System

This section is devoted to proving Theorem 2. The proof of this theorem is a generalization of the method of [8] and the idea is to produce a suitable supersolution and subsolution in the moving coordinate frame $z = x - ct$.

*Proof.* (Theorem 2) Let $(\psi(z), \phi(z), c)$ be as in Theorem 1 and choose $\bar{q}_u$, $\bar{q}_s$, $z_2$ such that

$$\liminf_{x \to -\infty} s_0(x) > 1 - \bar{q}_s > a \quad \liminf_{x \to -\infty} u_0(x) > g(1) - \bar{q}_u > g(a),$$

and

$$\psi(z + z_2) - \bar{q}_s \leq s_0(z) \quad \phi(z + z_2) - \bar{q}_u \leq u_0(z).$$

This is possible by choosing $z_2$ positive and large enough since $(s_0(z), u_0(z))$ satisfies (18). Define $\mathbf{v}_- := (w_-, v_-)^T$ where

$$w_- := \psi(z + \xi(t)) - q_s(t) \quad \text{and} \quad v_- = \phi(z + \xi(t)) - q_u(t).$$

The functions $\xi(t), q_s(t)$, and $q_u(t)$ will be defined later. Let $\tau = z + \xi(t)$ and apply the operator defined in (29) to $v_-$:

$$N\mathbf{v}_- = \begin{cases} \xi'(t)\psi'(\tau) - q_s'(t) - \psi''(\tau) - c\psi'(\tau) + \psi(\tau) - q_s(t) - \alpha(\phi(\tau) - q_u(t)) \\ \xi'(t)\phi'(\tau) - q_u'(t) - c\phi'(\tau) - g(\psi(\tau) - q_s(t)) + \phi(\tau) - q_u(t). \end{cases}$$

Using (17) we obtain

$$N\mathbf{v}_- = \begin{cases} \xi'(t)\psi'(\tau) - q_s'(t) - q_s(t) + \alpha q_u(t) \\ \xi'(t)\phi'(\tau) - q_u' + g(\psi(\tau)) - g(\psi(\tau) - q_s(t)) - q_u(t). \end{cases}$$

Now, choose $\delta > 0$ small enough such that for $\psi \in [0, \delta] \cup [1-\delta, 1]$ and $\phi \in [0, g(\delta)] \cup [g(1-\delta), g(1)]$ and $0 < q_s(t) \leq \bar{q}_s$ then

$$g(\psi(\tau)) - g(\psi(\tau) - q_s(t)) \leq \kappa q_s(t),$$

with $\kappa$ satisfying

$$\alpha \kappa < 1. \tag{31}$$

Such $\kappa$ exists given condition (14) as long as $\delta$ is chosen small enough. In the regime where $\xi'(t) > 0$ we obtain

$$N\mathbf{v}_- \leq \begin{cases} -q_s'(t) - q_s(t) + \alpha q_u(t) \\ -q_u'(t) + \kappa q_s(t) - q_u(t). \end{cases}$$

Then solving the system

$$q_s'(t) = -q_s(t) + \alpha q_u(t) \tag{32a}$$
$$q_u'(t) = -q_u(t) + \kappa q_s(t), \tag{32b}$$

with initial conditions $q_u(0) = \bar{q}_u$ and $q_s(0) = \bar{q}_s$ ($\bar{q}_s, \bar{q}_u$ defined above) we obtain solutions

$$q_s(t) = Ae^{\lambda_- t} + Be^{\lambda_+ t}$$
$$q_u(t) = \frac{-\sqrt{\alpha\kappa}}{\alpha}Ae^{\lambda_- t} + \frac{\sqrt{\alpha\kappa}}{\alpha}Be^{\lambda_+ t},$$

where $\lambda_\pm = -1 \pm \sqrt{\alpha\kappa} < 0$ by (31) for some constants $A$ and $B$. This guarantees that $q_s(t) \leq \bar{q}_s$ and $q_u(t) \leq \bar{q}_u$ and

$$\lim_{t \to \infty} q_s(t) = 0 \quad \text{and} \quad \lim_{t \to \infty} q_u(t) = 0.$$

Furthermore, (32) guarantees that $q_s(t), q_u(t) \geq 0$. We are left to treat the case when $(\psi, \phi) \in (\delta, 1-\delta) \times (g(\delta), g(1-\delta))$. Here $\psi', \phi' \leq -\gamma$ for some $\gamma > 0$. Since $\beta$ is the maximum value of $g'(s)$, our goal is to find a suitable $\xi(t)$ such that

$$N\mathbf{v}_- \leq \begin{cases} -\gamma \xi'(t) - q_s'(t) - q_s(t) + \alpha q_u(t) \\ -\gamma \xi'(t) - q_u'(t) + \beta q_s(t) - q_u(t). \end{cases}$$

Given our choice of $q_s(t)$ and $q_u(t)$ the first inequality above is automatically satisfied; therefore, we only need to solve the ODE

$$\gamma \xi'(t) = -q_u'(t) + \beta q_s(t) - q_u(t)$$
$$\overset{(32a)}{=} (\beta - \kappa) q_s(t),$$

with initial condition $\xi(0) = z_2$. The solution is

$$\xi(t) = z_1 + A_1 e^{\lambda_- t} + B_1 e^{\lambda_+ t},$$

where $A_1 = \frac{(\beta-\kappa)A}{\gamma \lambda_-}$ and $B_1 = \frac{(\beta-\kappa)B}{\gamma \lambda_+}$ and $z_1 = z_2 - (A_1 + B_1)$. Note that $\xi(t)$ is increasing and reaches a limit as $t$ approaches infinity; indeed,

$$\lim_{t \to \infty} \xi(t) = z_1.$$

Therefore, $z_2 \leq \xi(t) \leq z_1$ and at $t = 0$ we have

$$\psi(z + z_1) - \bar{q}_s \leq s_0(x) \quad \text{and} \quad \phi(z + z_1) - \bar{q}_u \leq u_0(x),$$

which guarantees that $\mathbf{v}_-$ is a subsolution. Also, since zero is a sub-solution then $\max\{0, \mathbf{v}_-\}$ is also a subsolution.

Similarly, we can check that $\phi(z + z_2) + q_u(t) \geq v(z,t)$ and $\psi(z + z_2) + q_s(t) \geq w(z,t)$. Thus, we obtain the upper and lower bounds in (19). Furthermore, as the upper and lower bounds for $u(x,t)$ in (19) approach $g(1)$ as $x, t \to \infty$. Finally, to prove (21) we take $\bar{q}_s, \bar{q}_u, |z_2 - z_1|$ of order $\mathcal{O}(\epsilon)$, which is possible by (20) and we conclude. $\square$

## 3.3 Unique Entire Solutions for the Heterogeneous Problem

This section is devoted to the proof of Theorem 3, which is our first result for the heterogenous problem ($L > 0$). A solution to (8), $\mathbf{v} = (w, v)^T$, satisfies $\mathcal{M}\mathbf{v} = 0$, where

$$\mathcal{M}\mathbf{v} = \begin{cases} w_t - w'' + w - \alpha_L(x) v \\ v_t - g(w) + v. \end{cases} \tag{33}$$

The proof of Theorem 3 requires some asymptotic estimates on the traveling wave solutions $(\psi, \phi)$, which we state below.

**Lemma 2** (Asymptotic Estimates). *Let $(\psi, \phi, c)$ be defined as in Theorem 1 then the following estimates hold for $\psi$ and $\phi$:*

$$\gamma e^{-\lambda z} \leq \psi(z), \phi(z) \leq \delta e^{-\lambda z} \quad z \geq 0 \tag{34a}$$
$$\gamma e^{\mu z} \leq 1 - \psi(z) \leq \delta e^{\mu z} \quad z < 0 \tag{34b}$$
$$\gamma e^{\mu z} \leq g(1) - g(\phi(z)) \leq \delta e^{\mu z}, \tag{34c}$$

*and for $\psi'(z)$ and $\phi'(z)$:*

$$-\tilde{\gamma} e^{-\lambda z} \leq \psi'(z), \phi'(z) \leq -\tilde{\delta} e^{-\lambda z} \quad z \geq 0 \tag{35a}$$
$$-\tilde{\gamma} e^{\mu z} \leq \psi'(z), \phi'(z) \leq -\tilde{\delta} e^{\mu z} \quad z < 0, \tag{35b}$$

*for $\gamma, \tilde{\gamma}, \tilde{\delta}, \mu, \lambda$ positive constants. Furthermore, $\lambda > \mu$.*

The proof of Lemma (2) is standard therefore we only prove $\lambda > \mu$ in the appendix. (see [8, 11] for the derivation of the estimates for the single equation). We mention, however, that the conclusion that $\lambda > \mu$, which is not automatic form standard exponential decay estimates, will be very important later.

As the traveling wave solutions $(\psi(z), \phi(z))$ are unique modulo translations, without loss of generality we set

$$\psi(0) = \frac{s_b}{2} \quad \text{and} \quad \phi(0) = \frac{g(s_b)}{2}. \tag{36}$$

The proof of Theorem 3 follow techniques of [4]. To prove existence it is useful to define, following [4], the auxiliary function $\xi(t)$ by

$$\xi(t) = \frac{1}{\lambda} \log \frac{1}{1 - c^{-1} M e^{\lambda ct}}, \tag{37}$$

where $c$ is the speed of the waves from Theorem 1, $\lambda$ is defined in Lemma 2, and $M$ a positive constant to be chosen later. Note that $\xi(t)$ is well-defined on $t \in (-\infty, -T)$ with $T := \frac{1}{\lambda c} \ln\left(\frac{c}{M+c}\right)$. Furthermore, $\lim_{t \to -\infty} \xi(t) = 0$ and

$$\dot{\xi}(t) = M e^{\lambda(ct+\xi)}. \tag{38}$$

For uniqueness we need to consider the spatial region where the front of the waves are located at a given time. More precisely, given $\eta \in [0, \frac{1}{2})$ we define the time function

$$\mathcal{F}_\eta(t) := \{x < 0 : \eta \leq \overline{s}(x,t) \leq 1 - \eta,\ g(\eta) \leq \overline{u}(x,t) \leq g(1-\eta)\}.$$

Since the waves are moving to the right there exists a time, $T_\eta \in \mathbb{R}$, such that $\mathcal{F}_\eta(t) \subset \{x \leq -1\}$ for $t \in (-\infty, -T_\eta)$. We state an extension of Lemma 3.1 in [4] and leave out the proof as it follows similarly to that of Lemma 3.1.

**Lemma 3.** *For any $\eta \in [0, 1/2)$, there exists a $\delta_\eta > 0$ such that for $(\overline{s}, \overline{u})$, the entire solutions to (8),*

$$\overline{s}_t, \overline{u}_t \geq \delta_\eta \quad \text{for } x \in \mathcal{F}_\eta(t),\ t \in (-\infty, T_\eta). \tag{39}$$

We now have all the tools to prove Theorem 3.

*Proof.* We break the proof into four steps. The first step is to show that the piecewise function, $\mathbf{v}_+ = (w_+, v_+)$, defined by

$$w_+(x,t) = \begin{cases} \psi(z_+) + \psi(z_-) & x < 0 \\ 2\psi(-ct - \xi(t)) & x \geq 0 \end{cases} \text{ and } v_+(x,t) = \begin{cases} \phi(z_+) + \phi(z_-) & x < 0 \\ 2\phi(-ct - \xi(t)) & x \geq 0, \end{cases} \tag{40}$$

where $z_+ = x - ct - \xi(t)$ and $z_- = -x - ct - \xi(t)$, with $\xi(t)$ satisfying (37), is a supersolution for $t \in (-\infty, -T)$, for $T$ large enough. Note that $w_+(x,t)$ and $v_+(x,t)$ are independent of the spatial variable in the gap and are both composed of a pair of self-annihilating fronts.

The second step is to find a suitable time range where $\mathbf{v}_- = (w_-, v_-)$ defined by

$$w_-(x,t) = \begin{cases} \psi(y_+) - \psi(y_-) & x \leq 0 \\ 0 & x > 0, \end{cases} \text{ and } v_-(x,t) = \begin{cases} \phi(y_+) - \phi(y_-) & x \leq 0 \\ 0 & x > 0, \end{cases} \tag{41}$$

for $y_+ = x - ct + \xi(t)$ and $y_- = -x - ct + \xi(t)$, is a subsolution for $t$ negative enough.

The third step is to take $\mathbf{v}_-$ at a suitably negative time as initial data and construct a monotone sequence of solutions, which are bounded below by the the subsolution and above by the supersolution for values of $t$ that are negative enough. From this we conclude the existence of an entire solution, that is bounded above by $\mathbf{v}_+$ and below by $\mathbf{v}_-$ for negative enough values of $t$. By the definition of the super and subsolutions then we obtain (22).

The final step is to show uniqueness and the proof will be a combination of the uniqueness proof found in [4] for the entire solution and the proof for Theorem 2.

*Step 1 (Supersolution):* We aim to show that $\mathbf{v}_+$ is a supersolution in a suitable time range.

*Case 1:* Consider the case when $x \leq 0$ and apply $\mathcal{M}$ (defined by (33)) to $\mathbf{v}_+$

$$\mathcal{M}\mathbf{v}_+ = \begin{cases} -(c+\xi'(t))(\psi'(z_+) + \psi'(z_-)) - (\psi''(z_+) + \psi''(z_-)) + \psi(z_+) + \psi(z_-) - \alpha(\phi(z_+) + \phi(z_-)) \\ -(c+\xi'(t))(\phi'(z_+) + \phi'(z_-)) - g(\psi(z_+) + \psi(z_-)) + \phi(z_+) + \phi(z_-) \end{cases}$$

$$\stackrel{(17)}{=} \begin{cases} -\xi'(t)(\psi'(z_+) + \psi'(z_-)) \\ -\xi'(t)(\phi'(z_+) + \phi'(z_-)) + g(\psi(z_+)) + g(\psi(z_-)) - g(\psi(z_+) + \psi(z_-)). \end{cases}$$

Since $\psi'(z)$ is always negative the first term in the above equality is strictly positive. For the second term we make use of the inequality (see [11])

$$|g(a) + g(b) - g(a+b)| \leq Lab, \tag{42}$$

for any $0 \leq a, b \leq 1$ and some constant $L > 0$, obtaining

$$\mathcal{M}\mathbf{v}_+ \geq \begin{cases} -\xi'(t)(\psi'(z_+) + \psi'(z_-)) \\ -\xi'(t)(\phi'(z_+) + \phi'(z_-)) - L\psi(z_+)\psi(z_-). \end{cases}$$

If $z_+ \leq 0$ then $|x| \geq -ct - \xi(t)$ as we are considering $x \leq 0$; thus, $z_- = -x - ct - \xi(t) \geq 0$. Thus, in this case $z_+ \leq 0$ and $z_- \geq 0$. We now invoke Lemma 2 we obtain

$$-\xi'(t)(\phi'(z_+) + \phi'(z_-)) - L\psi(z_+)\psi(z_-) \geq M\tilde{\delta}e^{\lambda(ct+\xi)}e^{\mu(x-ct-\xi)} - L\delta e^{-\lambda(-x-ct-\xi)}$$

$$\geq e^{\lambda(ct+\xi)}\left(M\tilde{\delta}e^{\mu z_+} - L\delta e^{\lambda x}\right).$$

Since $\lambda > \mu$ the above quantity is positive provided

$$M\tilde{\delta} > L\delta. \tag{43}$$

On the other hand, when $z_+ \geq 0$ (this implies that $z_- \geq 0$)

$$-\xi'(t)(\phi'(z_+) + \phi'(z_-)) - L\psi(z_+)\psi(z_-) \geq M\tilde{\delta}e^{\lambda(ct+\xi)}e^{-\lambda(x-ct-\xi)} - L\delta^2 e^{-\lambda(x-ct-\xi)}e^{-\lambda(-x-ct-\xi)}$$

$$\geq e^{2\lambda(ct+\xi)}(M\tilde{\delta}e^{-\lambda x} - L\delta^2).$$

Hence, for $x < 0$ if

$$M\tilde{\delta} > L\delta^2, \tag{44}$$

we have that $M\mathbf{v}_+ > 0$.

*Case 2:* If $x \geq 0$ then

$$\mathcal{M}\mathbf{v}_+ = \begin{cases} -2(c+\xi'(t))\psi'(-ct-\xi(t)) + 2\psi(-ct-\xi(t)) - 2\alpha_L(x)\phi(-ct-\xi(t)) \\ -2(c+\xi'(t))\phi'(-ct-\xi(t)) - g(2\psi(-ct-\xi(t))) + 2\phi(-ct-\xi(t)) \end{cases}$$

$$\stackrel{(17)}{=} \begin{cases} -2\xi'(t)\psi'(-ct-\xi(t)) \\ -2\xi'(t)\phi'(-ct-\xi(t)) + 2g(\psi(-ct-\xi(t))) - g(2\psi(-ct-\xi(t))). \end{cases}$$

Let $T_1 > 0$ be sufficiently large so that $t \in (-\infty, -T_1)$ then $-ct - \xi(t) \geq 0$, then $\psi(-ct - \xi) \leq \frac{s_b}{2}$. Since, $g(2\psi(-ct - \xi(t))) = 0$ in this case we have that $\mathbf{v}_+$ is a supersolution for $t \in (-\infty, -T_1)$ as $g(s) \leq 0$ for all $s \leq s_b$.

*Step 2 (Subsolution):* The case when $x > 0$ is trivial, hence, we assume that $x \leq 0$. Applying $\mathcal{M}$ to $\mathbf{v}_-$ we obtain

$$\mathcal{M}\mathbf{v}_- = \begin{cases} (-c+\xi')(\psi'(y_+) - \psi'(y_-)) - \psi''(y_+) + \psi''(y_-) + \psi(y_+) - \psi(y_-) - \alpha(\phi(y_+) - \phi(y_-)) \\ (-c+\xi')(\phi'(y_+) - \phi'(y_-)) - g(\psi(y_+) - \psi(y_-)) + \phi(y_+) - \phi(y_-) \end{cases}$$

$$\stackrel{(17)}{=} \begin{cases} \xi'(\psi'(y_+) - \psi'(y_-)) \\ \xi'(\phi'(y_+) - \phi'(y_-)) + g(\psi(y_+)) - g(\psi(y_-)) - g(\psi(y_+) - \psi(y_-)). \end{cases}$$

Observe that $y_+ \leq y_-$ as $x < 0$ and $y_- \geq 0$ because $x \leq 0$. There are two possibilities we need to consider, $y_+ \leq 0$ and $y_+ \geq 0$. In the latter case, both $\psi(y_+)$ and $\psi(y_-)$ are bounded above by $\frac{s_b}{2}$. First, this implies that $\psi'(y_+) - \psi'(y_-) \leq 0$ and $\phi'(y_+) - \phi'(y_-) \leq 0$. Second, it implies that $g(\psi(y_+)) - g(\psi(y_-)) - g(\psi(y_+) - \psi(y_-)) = 0$; hence, in this case $\mathcal{M}\mathbf{v}_- \leq 0$. We are left to consider the case when $y_+ \leq 0$, which we work out next (once again using (42))

$$\mathcal{M}\mathbf{v}_- \leq \begin{cases} Me^{\lambda(ct+\xi)} \left(-\tilde{\delta}e^{\mu y_+} + \tilde{\gamma}e^{-\lambda y_-}\right) := I \\ Me^{\lambda(ct+\xi)}(-\tilde{\delta}e^{\mu y_+} + \tilde{\gamma}e^{-\lambda y_-}) + L\psi(y_-)(\psi(y_+) - \psi(y_-)) := II. \end{cases}$$

Rewrite the first term

$$I = -Me^{\lambda(x+ct+\xi)} \left(\tilde{\delta}e^{\mu(-ct+\xi) + x(\mu-\lambda)} - \tilde{\gamma}e^{-\lambda(-ct+\xi)}\right)$$

$$\leq -Me^{\lambda(x+ct+\xi)} \left(\tilde{\delta}e^{\mu(-ct+\xi)} - \tilde{\gamma}e^{-\lambda(-ct+\xi)}\right).$$

Therefore, if $t < -T_1$, replace $T_1$ by a larger value if necessary, then

$$\tilde{\delta}e^{\mu(-ct+\xi(t))} - \tilde{\gamma}e^{-\lambda(-ct+\xi(t))} \geq 0 \quad \text{for } t \in (-\infty, -T_1] \tag{45}$$

because $\lambda > \mu$. For the second term we use the bound

$$L\psi(y_-)(\psi(y_+) - \psi(y_-)) \leq L\delta e^{-\lambda y_-},$$

which gives,

$$II \leq -Me^{\lambda(x+ct+\xi)} \left(\tilde{\delta}e^{\mu(-ct+\xi)} - \tilde{\gamma}e^{-\lambda(-ct+\xi)} - \frac{L}{M}\delta e^{-2\lambda\xi}\right)$$

$$\leq 0,$$

provided $T_1$ is chosen negative enough to satisfy

$$\tilde{\delta}e^{\mu(-ct+\xi)} - \tilde{\gamma}e^{-\lambda(-ct+\xi)} - \frac{L}{M}\delta e^{-2\lambda\xi} \geq 0 \quad \text{for } t \in (-\infty, -T_1],$$

with $T_1$ replaced by its updated value.

*Step 3 (Entire Solution):* We obtain the desired entire solutions by constructing a sequence of solutions that will be defined on the time interval $(-n, \infty)$. As $n \to \infty$ these sequences converge uniformly in $x$ as $n \to \infty$. Let $n$ be a sufficiently large positive integer and $(s_n(x,t), u_n(x,t))$ be the solutions to (8) for $x \in (-\infty, \infty)$ and $t \in (-n, \infty)$ with initial data

$$s_n(x, -n) = w_-(x, -n), \quad u_n(x, -n) = v_-(x, -n).$$

From the definitions of $\mathbf{v}_+, \mathbf{v}_-$, see (40) and (41) respectively, we see that

$$w_-(x, -n) \leq s_n(x, -n) \leq w_+(x, -n) \quad \text{and} \quad v_-(x, -n) \leq v_n(x, -n) \leq v_+(x, -n)$$

By the comparison principle

$$w_-(x,t) \leq s_n(x,t) \leq w_+(x,t), \ v_-(x,t) \leq u_n(x,t) \leq v_+(x,t) \text{ for } (x,t) \in \mathbb{R} \times (-n, -T_1), \quad (46)$$

with $-T_1$ defined at the end of *Step 2*. In particular, (46) holds for $t = -n + 1 = -(n-1)$, which implies that for $x \in \mathbb{R}$

$$s_{n-1}(x, -n+1) = w_-(x, -n+1) \leq s_n(x, -n+1),$$
$$u_{n-1}(x, -n+1) = v_-(x, -n+1) \leq u_n(x, -n+1).$$

Another application of the comparison principle gives that

$$s_{n-1}(x,t) \leq s_n(x,t) \text{ and } u_{n-1}(x,t) \leq u_n(x,t) \text{ for } (x,t) \in \mathbb{R} \times (-n+1, -T_1).$$

Thus, we obtain monotone increasing sequences $\{s_n\}_{n \in \mathbb{Z}^+}$ and $\{u_n\}_{n \in \mathbb{Z}^+}$, which enable us take the limit as $n \to \infty$, uniformly in $x$, to obtain solutions $(\overline{s}(x,t), \overline{u}(x,t))$ to (8) on $(x,t) \in \mathbb{R} \times \mathbb{R}$. Furthermore, we have the additional bounds

$$w_-(x,t) \leq \overline{s}(x,t) \leq w_+(x,t), \ v_-(x,t) \leq \overline{u}(x,t) \leq v_+(x,t) \text{ for } (x,t) \in \mathbb{R} \times (-\infty, -T_1).$$

The definitions of $\mathbf{v}_+, \mathbf{v}_-$ also give (22). We are left to prove that

$$\overline{s}_t(x,t), \ \overline{u}_t(x,t) > 0. \quad (47)$$

We begin by looking at the time derivative of $\mathbf{v}_-$ when $x < 0$

$$\partial_t w_- = (-c + \xi'(t))(\psi'(y_+) - \psi'(y_-)) \quad \partial_t v_- = (-c + \xi'(t))(\phi'(y_+) - \phi'(y_-)).$$

Take $t$ negative enough so that $-c + \xi'(t) < 0$ and arguing as before $(\psi'(y_+) - \psi'(y_-)) < 0$. Note that if $y_+ \geq 0$ this inequality is clear and if $y_+ \leq 0$ then we take $t$ negative enough. Thus, $\partial_t w_- > 0$ for $t$ sufficiently negative (a similar argument works for $v_-$), which by (46) implies that $\partial_t s_n(x,t), \partial_t u_n(x,t) > 0$ if $n$ is sufficiently large as $u_n(x,-n) = v_-(x,-n)$ (similarly for $s_n$ and $w_-$). Applying the maximum principle gives that

$$\partial_t s_n(x,t) > 0, \ \partial_t u_n(x,t) > 0 \text{ for } t \in (-n, \infty).$$

Taking the limit as $n \to \infty$ and using the fact that $\partial_t \overline{s}(x, -\infty) > 0$ (same for $\partial_t \overline{u}$) proves (47).

*Step 4 (Uniqueness):* Assume for contradiction that $(s(x,t), u(x,t))$ is another solution to (8) that also satisfies (22). Since $(\overline{s}, \overline{u})$ and $(s, u)$ both satisfy (22) then for any $\epsilon > 0$ there exists a time $t_\epsilon \in \mathbb{R}$ such that

$$\|\overline{s}(\cdot, t) - s(\cdot, t)\|_\infty < c\epsilon \quad \|\overline{u}(\cdot, t) - u(\cdot, t)\|_\infty < c\epsilon \quad (48)$$

for all $t \in (-\infty, t_\epsilon)$, where $c = \min\{\overline{q}_s, \overline{q}_u\}$ ($\overline{q}_s, \overline{q}_u$ defined as in the proof of Theorem (2)). For $t_0 \in (-\infty, t_\epsilon - \epsilon]$ define

$$S^+(x,t) := \overline{s}(x, t + t_0 + \epsilon\xi(t)) + \epsilon q_s(t), \quad U^+(x,t) := \overline{u}(x, t + t_0 + \epsilon\xi(t)) + \epsilon q_u(t)$$
$$S^-(x,t) := \overline{s}(x, t + t_0 - \epsilon\xi(t)) - \epsilon q_s(t), \quad U^-(x,t) := \overline{u}(x, t + t_0 - \epsilon\xi(t)) - \epsilon q_u(t).$$

As in the proof of Theorem 2 we aim find $\xi(t)$, $q_u(t)$ and $q_s(t)$ such that $(S^+, U^+)$ is a supersolution and $(S^-, U^-)$ is a subsolution for $t \in [0, T_\eta - t_0 - \epsilon]$ with $T_\eta$ defined in (39) for $\eta$ to be fixed later. The proof is very similar to that of Theorem 2 and we only summarize the steps for the subsolution and note that the supersolution is checked similarly. Let $\tau = t + t_0 - \epsilon\xi(t)$ and apply the operator $\mathcal{M}$ (defined in (33)) to $(S^-, U^-)$:

$$\partial_t S^- - S^-_{xx} + S^- - \alpha U^- = -\epsilon(\xi'(t)\partial_t \overline{s}(x,\tau) + q'_s(t) + q_s(t) - \alpha q_u(t))$$
$$\partial_t U^- - g(S^-) + U^- = -\epsilon\xi'(t)\partial_t \overline{u}(x,\tau) - \epsilon q'_u(t) + g(\overline{s}(x,\tau)) - g(\overline{s}(x,\tau) - \epsilon q_s(t)) - \epsilon q_u(t).$$

In the region outside the front of the wave, $x \notin \mathcal{F}_\eta(t + t_0 - \epsilon\xi(t))$, then for $\eta$ small enough we have that $g(\bar{s}(x,\tau)) - g(\bar{s}(x,\tau) - \epsilon q_s(t)) \leq \kappa q_s(t)$, where once again we choose $\kappa$ such that (31) holds. Then, as in the proof of Theorem 2 we can find positive and bounded functions $q_u(t)$ and $q_s(t)$ that satisfy (32). Recall, that both $q_s(t), q_u(t) \to 0$ as $t \to \infty$. On the other hand, when $x \in \mathcal{F}_\eta(t + t_0 - \epsilon\xi(t))$ we invoke (39), noting that this holds when $t \leq T_\eta - t_0 - \epsilon$, and solve for $\xi(t)$

$$\xi'(t) = \frac{\beta - \kappa}{\delta_\eta} q_s(t) \quad \xi(0) = 0,$$

with $\delta_\eta$ satisfying (39). Here, since $\beta > \kappa$ we also get that $\xi(t)$ increases but approaches a limit as $t \to \infty$, say

$$\lim_{t \to \infty} \xi(t) = M.$$

Thus, $(S^-, U^-)$ is subsolution and $(S^+, U^+)$ is a supersolution for $t \in [0, T_\eta - t_0 - M\epsilon]$. Inequality (48) implies that

$$S^-(x, 0) \leq s(x, t_0) \leq S^+(x, 0) \quad \text{and} \quad U^-(x, 0) \leq u(x, t_0) \leq U^+(x, 0).$$

Implying that for $t \in [0, T_\eta - t_0 - \epsilon]$ then

$$S^-(x, t) \leq s(x, t + t_0) \leq S^+(x, t) \text{ and } U^-(x, t) \leq u(x, t + t_0) \leq U^+(x, t)$$

for all $x \in \mathbb{R}$. By a change of variables $\tilde{t} = t + t_0$, dropping the tilde for simplicity, this is equivalent to

$$\bar{s}(x, t - \epsilon\xi(t - t_0)) - \epsilon q_s(t - t_0) \leq s(x, t) \leq \bar{s}(x, t + \epsilon\xi(t - t_0)) + \epsilon q_s(t - t_0)$$
$$\bar{u}(x, t - \epsilon\xi(t - t_0)) - \epsilon q_u(t - t_0) \leq u(x, t) \leq \bar{u}(x, t + \epsilon\xi(t - t_0)) + \epsilon q_u(t - t_0),$$

for $t \in [t_0, T_\eta - M\epsilon]$ and $t_0 \in (-\infty, T_\eta - M\epsilon]$. Taking the limit $t_0 \to -\infty$, using the facts that $q_u(t - t_0), q_s(t - t_0) \to 0$ and $\xi(t) \to M$, we obtain

$$\bar{s}(x, t - M\epsilon) \leq s(x, t) \leq \bar{s}(x, t + M\epsilon)$$
$$\bar{u}(x, t - M\epsilon) \leq u(x, t) \leq \bar{u}(x, t + M\epsilon)$$

for $(x, t) \in \mathbb{R} \times (-\infty, T_\eta - M\epsilon]$. Finally, letting $\epsilon \to 0$ gives that $s \equiv \bar{s}$ and $u \equiv \bar{u}$ and we conclude. $\square$

### 3.4 Preventing the Propagation of Criminal Activity

This section is devoted to the problem of preventing the invasion of the criminal activity. Consider the alternative representation of the steady state solution. Any solution to (24) satisfies, for any $y_1, y_2 \in (-\infty, 0]$ or $y_1, y_2 \in [L, \infty)$,

$$\frac{1}{2}\left[(s'_L(y_1))^2 - (s'_L(y_2))^2\right] = -\int_{y_1}^{y_2} s''_L(x) s'_L(x) \, dx = \int_{y_1}^{y_2} f(s_L(x)) s'_L(x) \, dx = \int_{s_L(y_1)}^{s_L(y_2)} f(\theta) \, d\theta.$$

Letting $y_1 = -\infty$ and $y_2 = x$ above then (23) implies

$$\frac{1}{2}(s'_L(x))^2 = \int_{s_L(x)}^{1} f(\theta) \, d\theta.$$

for any $x \leq 0$. On the other hand, when $x \in [0, L]$ the solution must have the form

$$s_L(x) = Ae^{-x} + Be^x$$

for some $A, B \in \mathbb{R}$. As we seek $C^1$ solutions the following must hold:

$$\begin{cases} \frac{1}{2}(B - A)^2 = \mathcal{F}(1) - \mathcal{F}(A + B) \\ s_L(x) = Ae^{-x} + Be^x \\ \frac{1}{2}(-Ae^{-L} + Be^L)^2 = \mathcal{F}(Ae^{-L} + Be^L) - \mathcal{F}(s_L(\infty)), \end{cases} \quad (49)$$

for $0 \leq x \leq L$. If there are values $L \in \mathbb{R}^+$ and $A, B \in \mathbb{R}$ with the bound $0 \leq A + B \leq 1$, then one can explicitly build a solution to (24) using the quadrature method (see [7]) for appropriate values of $s_L(\infty)$. Let $\psi_L(x)$ be the unique solution to

$$u'' + f(u) = 0 \quad x \in (-\infty, 0)$$
$$u(-\infty) = 1 \quad u(0) = A + B.$$

Such a solution exists as $A + B = u(0) \in [0, 1]$ and it is monotonically decreasing; in fact,

$$\psi_L'(x) = -\sqrt{2 \int_{\psi_L(x)}^1 f(\theta) \, d\theta}.$$

Thus, $\psi_L'(x) < 0$ because the right hand side of the above equation is bounded away from zero. Thus, if a unique solution, $\psi_R(x)$, exists to

$$u'' + f(u) = 0 \quad x \in (L, \infty) \tag{50a}$$
$$u(L) = Ae^{-L} + Be^L \quad u(\infty) = s_L(\infty), \tag{50b}$$

then we can build an explicit solution of the form

$$s_L(x) = \begin{cases} \psi_L(x) & x \in (-\infty, 0] \\ Ae^{-x} + Be^x & x \in (0, L) \\ \psi_R(x) & x \in (L, \infty]. \end{cases}$$

Note that (50) has solutions only for certain values of $s_L(\infty)$.

*Remark* 1. To find solutions to (24) it suffices to find appropriate values of $L, A, B$ and a unique $\psi_R(x)$ satisfying (50). A blocking solution will necessarily have that $s_L(\infty)$ is bounded away from one.

Now we move to the proof of Theorem 4. In order to prevent the propagation of the waves the existence of a steady state solution that is uniformly bounded away from one as $x \to \infty$ is necessary. In fact, we will show that with sufficient resources a strictly decreasing solution that approaches zero as $x \to \infty$ will always block the propagation of the wave. We first show that for $L$ sufficiently large such a steady state solution exist.

### 3.4.1 Case I: Sufficient Resources

Our goal in this subsection is to build super and subsolutions to (24) for $L$ large enough, where both the supersolution and subsolution decay to zero at positive infinity. In particular, we seek a supersolution, $\phi_+(x)$, such that $\phi_+(0) = 1$ and $\phi_+(L) = b$, where $b$ is defined by $\mathcal{F}(b) = 0$. To be precise, let $\phi_+(x) = Ae^{-x} + Be^x$ for $x \in [0, L]$, we seek to find $A$, $B$, and $L$ such that

$$A + B = 1, \ Ae^{-L} + Be^L = b, \ \text{and} \ -Ae^{-L} + Be^L = 0.$$

This reduces our problem to finding $L$ such that $\cosh L = \frac{1}{b}$, which always has a solution because $0 < b < 1$. The unique values are given by

$$L_0 = \cosh^{-1}\left(\frac{1}{b}\right) > 0, \; B = \frac{b}{2}e^{-L_0}, \text{ and } A = \frac{b}{2}e^{L_0}.$$

We use $L_0$ to construct a monotonically decreasing solution to (24) with $L = L_0$.

**Lemma 4** (Base Case)**.** *There exists a monotonically decreasing solution, $s_{L_0}(x)$, to (24) with $L = L_0$, such that*

$$\lim_{x \to -\infty} s_{L_0}(x) = 1 \quad \lim_{x \to \infty} s_{L_0}(x) = 0.$$

*Proof.* We prove this lemma in three steps. The first step is to construct a suitable supersolution and subsolution. We then build a solution by using the supersolution and subsolution as barriers. The first solution constructed is not necessarily monotone, but it is bounded above by $b$ for all $x > L$. Using this upper bound we prove that there also exists a monotone solution that is below the first solution we constructed.

*Step 1: (Super/subsolution)* Let $\psi(x)$ be the unique solution to Cauchy problem

$$-u'' = f(u) \quad x \in [L_0, \infty)$$
$$u(L_0) = b, \quad u'(L_0) = 0.$$

Such a solution exists and in fact following the arguments from the beginning of §3.4 we can show that $\psi'(x) < 0$ and $\lim_{x \to \infty} \psi(x) = 0$. Now, define $\phi_+(x)$ by

$$\phi_+(x) = \begin{cases} 1 & x \leq 0 \\ \frac{b}{2}\left(e^{L_0}e^{-x} + e^{-L_0}e^x\right) & 0 \leq x \leq L_0 \\ \psi(x) & x > L_0. \end{cases}$$

The function $\phi_+(x)$ is a supersolution as it satisfies

$$\lim_{x \to \xi^-} \phi'_+(x) \geq \lim_{x \to \xi^+} \phi'_+(x),$$

for $\xi = 0$ and $\xi = L_0$. Similarly, we build a subsolution, $\phi_-(x)$. Let $\psi_-$ now be the solution to the Cauchy problem

$$-u'' = f(u) \quad x \in [-\infty, 0)$$
$$u(-\infty) = 1, \; u(0) = 0.$$

As before, we know that such a solution exists, and we define the subsolution

$$\phi_- = \begin{cases} \psi_-(x) & x \leq 0 \\ 0 & x > 0. \end{cases}$$

Furthermore, $\phi'_-(0) = -\sqrt{2\int_0^1 f(\theta)d\,\theta}$, which is bounded away from zero and so the condition on the derivative at $x = 0$ for a subsolution is satisfied.

*Step 2: (Barrier Method)* Having constructed a supersolution and a subsolution we use the standard barrier method (see for example [18]) to show that there exists a solution to (24). This steps follow those of the proof of the existence of the entire solutions in Theorem 3 (see *Step 3* of the proof), so we leave out the details.

*Step 3:* If the solution is monotone then we are done. If it is not monotone then the lack of monotonicity happens in the gap. Indeed, for $x > L$ we know that $s_L(x) \leq b$ and then

$$-\frac{1}{2}(s'_L(x))^2 = \int_0^{s_L(x)} f(\theta) \, d\theta < 0.$$

Since, the solution must be concave up in the gap then it reaches a minimum at a unique point $x^\star \in (0, L)$. If $s_L(x^\star) \geq a$ we construct a supersolution by letting $s_+(x) = s_L(x)$ for $x \leq x^\star$, $s_+(x) = s_L(x^\star)$ for $x \in [x^\star, x_1]$, where $x_1 > L$ is the unique value such that $s_L(x_1) = s_L(x^\star)$. For $x > x_1$ then $s_+(x) = s_L(x)$. Hence, $s_+(x)$ cuts off the non-monotone part of the solution $s_L(x)$. The other possibility is that $s_L(x^\star) < a$ and in this case cut off from $[x, y]$ with $x < L$ and $y > L$ being the unique point with $s_L(x) = s_L(y) = a$. In both cases, $s_+(x)$ is a supersolution. Once more the barrier method will allow us to construct a solution $s_{L_0}(x)$, which is monotone. Indeed, in the latter case, $s_{L_0}(x) \leq a$ for $x \geq x^\star$ and $s'_{L_0}(x) < 0$ in this regime. In the former case, $s_{L_0}(x)$ cannot have a non-negative derivative everywhere because it can only be zero (and thus later negative) if $s_{L_0}(x) = b$, which can never happen. For $x \in (0, x^\star)$ the solution is monotone by definition. Furthermore, the requirement that in the limit $x \to -\infty$ the solution $s_L(x) \to 1$ guarantees monotonicity for $x < 0$ (see the proof of Lemma (7) part c). This concludes the proof. □

The maximum value of the solution at at $L$ is indeed $b$, which is the maximum value that the homoclinic orbit through the point $(s(x) = 0, s'(x) = 0)$ attains. For $L = L_0$ is is clear that the steady state solution has a value less than or equal to $b$ at $L_0$ as the solution must remain below the supersolution. This brings up the question of whether there exists an $L < L_0$ such that $s(L) = b$ and $s'(L) = 0$, so that the maximum value can actually be achieved? The answer to this question is provided in the following lemma.

**Lemma 5** (Necessary Resources). *There exists a unique $L^\star > 0$ such that there exists a monotone solution, $s_{L^\star}(x)$, to (24) with $s_{L^\star}(L^\star) = b$ and $s'_{L^\star}(L^\star) = 0$. Furthermore, any solution, $s_L(x)$, to (24), with $L < L^\star$, must satisfy $s_L(L) > s_{L^\star}(L^\star) = b$.*

*Proof.* The first objective is to find $L$ such that $s_L(L) = b$ and $s'_L(L) = 0$ so that we know that $Ae^{-L} + Be^L = b$ and $Ae^{-L} = Be^L$. Note that this immediately satisfies the last equality in (49) with $s_L(\infty) = 0$. Substituting this into the first equality of (49) gives

$$\frac{1}{2}(b \sinh L)^2 = \mathcal{F}(1) - \mathcal{F}(b \cosh L),$$

which has a solution. Indeed, for $L = 0$ the right hand side is positive and the left hand side is zero. On the other hand, for $L_0 = \cosh^{-1}(\frac{1}{b})$ the left hand side is positive and the right hand side is zero. Hence, by continuity there exists a solution, which we call $L^\star \in (0, L_0)$. Given $L^\star$ we have a solution

$$s_{L^\star}(x) = \frac{b}{2}(e^{L^\star - x} + e^{x - L^\star}) \text{ for } x \in [0, L^\star],$$

which is monotonically decreasing and reaches a minimum at $x = L^\star$. For notational simplicity we denote $s_{L^\star}$ by $s^\star$. For any $x > L^\star$ we have

$$\frac{1}{2}(\partial_x s^\star(x))^2 = \int_{s^\star(x)}^b f(\theta) \, d\theta,$$

which implies that $s^\star(x) < b$ for all $x > L^\star$. Furthermore, since $\partial_x s^\star(x) \neq 0$ when $s^\star \in (0, b)$, $s^\star(x)$ cannot be a constant on any interval, and $\partial_x s^\star(x)$ is continuous then $\partial_x s^\star(x) < 0$ for all $x > L^\star$ with $s^\star(x) \neq 0$. On the left side of the interval, $s^\star(0) > b$ and so for $x \in (-\infty, 0]$ $s^\star(x) \in (b, 1)$ and $(s^\star)'(x) < 0$. We are left to prove that

$$\lim_{x \to -\infty} s^\star(x) = 1 \quad \text{and} \quad \lim_{x \to +\infty} s^\star(x) = 0.$$

This is true, of course, only because $s^\star(x) > b$ for $x < 0$ and $s^\star < b$ for $x > L^\star$.

*Claim:* $\lim_{x \to -\infty} s^\star(x) = 1$. We know that the solution is bounded by one and it is monotonic. Thus, we only prove that limit is not somewhere in between $b$ and $1$. This is accomplished by proving that the reaction term $f(u)$ approaches zero. Consider the function

$$v(x) = \frac{1 - x^2}{2},$$

which is the solution to

$$-v'' = 1 \text{ in } (-1, 1) \quad v = 0 \text{ on } \{-1, 1\}.$$

For $x < 0$, define

$$\Gamma(x) = \inf \{f(\theta) : \theta \in [b, s_L^\star(x)]\}.$$

We will prove that $\Gamma(x) \leq \frac{2}{\text{dist}(x,0)^2}$, which then implies that $f(s)$ vanishes as $x \to -\infty$. Assume for contradiction that this is not the case, so that there exists a $x_0 \in (-\infty, 0)$ such that

$$\Gamma(x_0) > \frac{2}{\text{dist}(x_0, 0)^2}. \tag{51}$$

Note that $x_0$ cannot be a local minimum, as $\partial_{xx} s_L^\star(x_0) = -f(s_L^\star(x_0)) < 0$, thus we can find $x_1$ close to $x_0$ such that $s_L^\star(x_1) < s_L^\star(x_0)$ and $r < \text{dist}(x_1, 0) < \text{dist}(x_0, 0)$. Here, $r$ is chosen so that (51) still holds, that is,

$$\Gamma(x_0) > \frac{2}{r^2}.$$

Now consider the solution,

$$z(x) = \Gamma(x_0) r^2 v\left(\frac{x - x_1}{r}\right),$$

to

$$-z'' = \Gamma(x_0) \text{ for } x \in I_r(x_1)$$
$$z = 0 \text{ on } \partial I_r,$$

where $I_r(x_1)$ is the interval centered on $x_1$ with radius $r$. We now scale $z(x)$ by $\tau$ by defining

$$z^\tau = \tau z(x).$$

It is clear that for $\tau$ small enough $z^\tau(x) \leq s_{L^\star}(x)$ for $x \in I_r(x_1)$ and as $\tau$ increases there is a first $\tau^\star$ and $x^\star$ where $z^{\tau^\star}$ touches the graph of $s_L^\star(x)$, that is $s_{L^\star}(x^\star) = z^{\tau^\star}(x^\star)$. We know also that $x^\star$ is inside the interval. Thus, using the fact that the maximum of $z(x)$ is achieved at $x_1$, we obtain

$$s_{L^\star}(x^\star) = z^{\tau^\star}(x^\star) \leq \tau^\star z(x_1) = \frac{\tau^\star \Gamma(x_0) r^2}{2} \leq s_{L^\star}(x_1) < s_{L^\star}(x_0) < 1.$$

Thus, $\tau^\star < \frac{2}{\Gamma(x_0) r^2} < 1$, by our assumption on $r$. Define $w(x) := \tau^\star z(x) - s_L^\star(x)$ on $I_r(x_1)$. From the above argument we know that

$$w(x) \leq 0 \text{ for } x \in I_r(x_1) \text{ and } w(x^\star) = 0.$$

Now, since $s_{L^\star}(x^\star) < s_{L^\star}(x_0)$ there exists a neighborhood of $x^\star$, $\mathcal{V}$, such that $s_L^\star(x) < s_L^\star(x_0)$ for all $x \in \mathcal{V}$. This implies

$$s_{L^\star}''(x) \leq -\Gamma(x_0) \text{ for } x \in \mathcal{V}.$$

Thus, in this neighborhood, $\mathcal{V}$,

$$w''(x) \geq \Gamma(x_0)(1-\tau) > 0.$$

However, this contradicts the fact that $w(x)$ has a local maximum at $x^\star$. Thus, $\Gamma(x) \to 0$ which implies that $f(s) \to 0$ and from this we conclude the claim. The remaining limit, $\lim_{x \to +\infty} s_L^\star(x) = 0$ is proved similarly by looking at the equation that $1 - s_L^\star(x)$ solves.

Now, consider the case when $L < L^\star$ clearly $s^\star(x)$ is a subsolution to the problem

$$s'' + f_L(x, s) = 0. \tag{52}$$

Any solution to (52), $s_L(x)$, must satisfy $s_L(x) \geq s^\star(x)$; hence, it is true that $s_L(L) \geq s^\star(L) > s^\star(L^\star)$ and in particular $s_L(L) > b$. Similarly, for $L > L^\star$, $s^\star(x)$ is a supersolution to (24) and thus $s_L(L) \leq s^\star(L) < s^\star(L^\star)$. Hence, $L^\star$ is unique. □

Next, we show that for any $L$ larger than $L^\star$ there also exists a monotone decreasing solution.

**Lemma 6** (Sufficient Resources). *There exists a monotonically decreasing solution, $s_L(x)$, to (24) when $L > L^\star$, such that*

$$\lim_{x \to -\infty} s_L(x) = 1 \qquad \lim_{x \to \infty} s_L(x) = 0.$$

*Proof.* Having constructed a solution for $L^\star$, call it $s^\star(x)$ we now show that for $L > L^\star$ there also exists a steady state solution. Indeed, note that we have

$$s_{xx}^\star + f_{L^\star}(x, s^\star) = 0,$$

and so $s^\star$ is a supersolution to the problem,

$$s'' + f_L(x, s) = 0$$

as $f_L(x, s^\star) \leq f_{L^\star}(x, s^\star)$. The rest of the proof follows the proof of Lemma 4. □

### 3.4.2 Case II: Insufficient Resources

Now, we treat the case when $L$ is small and the wave propagates. Indeed, when $L = 0$ we know that for appropriate initial conditions the wave will propagate, see Theorem 1. In the first lemma we prove the existence of a symmetric solution, connecting one to itself, for any $L > 0$. In the case when $L = 0$ the symmetric solution is $s(x) \equiv 1$.

**Lemma 7** (Symmetric Steady State Solution). *For all $L > 0$ there exists a solution, $s_m(x)$, to (10b), with the following properties:*

(a) $s_m(x)$ is symmetric about $x = L/2$.

(b) $\lim_{x \to \pm\infty} s_m(x) = 1$.

(c) $s_m(x)$ is monotone increasing on $[L, \infty)$.

*Proof.* Since this problem is translation invariant, we prove the existence of a symmetric solution by translating the problem to the left by $L/2$. In other words,

$$\alpha(x) = \begin{cases} \alpha & x \in (-\infty, \frac{-L}{2}] \cap [\frac{L}{2}, \infty) \\ 0 & x \in (\frac{-L}{2}, \frac{L}{2}), \end{cases}$$

so that any symmetric solution must be of the form $s_m(x) = Ae^{-x} + Ae^x$ on $[\frac{-L}{2}, \frac{L}{2}]$. As noted before the steady state problem is equivalent to (49) with $\mathcal{F}(s_L(\infty)) = \mathcal{F}(1)$, and so the symmetric condition implies that the following equality must hold

$$\frac{1}{2}\left(2A\sinh\left(\frac{L}{2}\right)\right)^2 = \mathcal{F}(1) - \mathcal{F}\left(A\cosh\left(\frac{L}{2}\right)\right).$$

Letting $z = 2A\cosh(\frac{L}{2})$ the above reduces to

$$z^2 \tanh^2\left(\frac{L}{2}\right) = \mathcal{F}(1) - \mathcal{F}(z),$$

which always has a solution as $\mathcal{F}(1) > 0$. Hence, we can build explicitly the symmetric solution about the origin by the method discussed earlier (see Remark 1) such that

$$\lim_{x \to \pm\infty} s_m(x) = 1.$$

We are left to show that this solution is monotone outside of $[-L/2, L/2]$. To do this we prove that $s'_m(x) > 0$ for $x > L/2$ given the fact that $\lim_{x \to \infty} s_m(x) = 1$ (by symmetry then $s'_m(x) < 0$ for $x < -L/2$). Choose $\delta$ such that $s_m(L/2) < 1 - \delta$, then the limit implies that for $\delta > 0$ there is a constant $R > 0$, sufficiently large, such that

$$s_m(x) > 1 - \delta \text{ for } x > R.$$

Consider a translated solution

$$s_m^\tau(x) = s_m(x + \tau) \text{ for } \tau > 0.$$

Note that $s^\tau(x) > 1 - \delta$ for $x > R - \tau$, so for $\tau = R - L/2$ we obtain

$$s_m^\tau(L/2) > s_m(L/2).$$

Invoking the maximum principle for semi-infinite domains for bistable reaction-diffusion systems, see [2], we obtain that

$$s_m^\tau(x) > s_m(x) \quad \text{for} \quad x > L/2.$$

Now, define

$$\tau^\star = \inf\left\{\tau > 0 \mid s_m^\tau(x) \geq s_m(x) \text{ for } x > L/2\right\},$$

it is clear that $\tau^\star < R - L/2$. Let us assume for contradiction that $\tau^\star > 0$, by definition of $\tau^\star$ then

$$G := \inf_{I_a}(s_m^{\tau^\star}(x) - s_m(x)) \geq 0, \quad \text{with } I_a = [L/2, R].$$

If $G > 0$ then, by standard elliptic theory [10], there exists a $\epsilon > 0$ such that

$$s_m^t(x) > s_m(x) \quad x \in I_a, \ t \in (\tau^\star - \epsilon, \tau^\star).$$

The maximum principle for unbounded domains on $[R, \infty)$ gives that for some $\epsilon_0 \in (0, \epsilon)$,

$$s_m^{\tau^\star - \epsilon_0}(x) > s_m(x) \quad x \in [L/2, \infty).$$

This contradicts the definition of $\tau^\star$. If $G = 0$, there exists a sequence $\{\xi_k\}_{k \in \mathbb{N}} \in I_a$ such that

$$\lim_{k \to 0}\left(s_m^{\tau^\star}(\xi_k) - s_m(\xi_k)\right) = 0.$$

We translate the solutions by defining $s_m^k(x) = s_m(x + \xi_k)$, and obtain a bounded sequence of solutions. Thus, a subsequence converges to a limit, $\bar{s}_m(x)$, on compact sets. Furthermore, as $G = 0$ we have that $\bar{s}_m^{\tau^\star}(0) = \bar{s}_m(0)$ and $s_m^k(x) \geq s_m(x)$ for all $k$; thus, $\bar{s}_m(x + \tau^\star) = \bar{s}_m(x)$ for all $x > 0$ by the strong maximum principle. Thus, $\bar{s}_m(x)$ must be a periodic function with period $\tau^\star$. However, recall that $\xi_k \in [L/2, R]$ and so $\bar{s}_m(x) \to 1$ as $x \to \infty$, which gives a contradiction. Thus, $\tau^\star = 0$ and we conclude that $s_m(x + t) > s_m(x)$ for all $t \geq 0$.

□

We now show that for small $L$ there are no solutions that are bounded away from one for $x > L$.

**Lemma 8.** *For $L < L^\star$ there are no non-negative solutions, $s_L(x)$, to (10b) that are bounded above by one, and such that*

$$\lim_{x \to -\infty} s_L(x) = 1 \quad \lim_{x \to \infty} s_L(x) \leq 1 - \epsilon. \tag{53}$$

*for any $\epsilon > 0$ small.*

*Proof.* Assume, for contradiction, that there exists a solution $s_L(x)$ to (10b) such that it satisfies (53). Since $L < L^\star$ from Lemma 5 we know that $s_L(L) > b$. Furthermore, for $L < y < x$

$$\left(s'_L(y)\right)^2 - \left(s'_L(x)\right)^2 = 2 \int_{s_L(y)}^{s_L(x)} f(\theta) \, d\theta. \tag{54}$$

Note that the integral is bounded away from zero as $s_L(x) \in (b, 1)$; thus, $s'_L(x) \neq s'_L(y)$ and the derivative is monotone. Assume for contradiction that $s'_L(x) < 0$ for all $x > L$, then $s_L(x) < s_L(L)$, from (54) we get that $s'_L(L) \leq s'_L(x)$. Thus, $\lim_{x \to \infty} s'_L(x) = 0$, if $s_L(x)$ is to stay above $b$. Looking at $x \to \infty$, we obtain

$$\left(s'_L(y)\right)^2 = 2 \lim_{x \to \infty} \int_{s_L(y)}^{s_L(x)} f(\theta) \, d\theta < 0.$$

This is a contradiction, and so it must be that $s'_L(x) > 0$ and by prior arguments, see the proof of Theorem 5, we know that

$$\lim_{x \to \infty} s_L(x) = 1.$$

This gives the final contradiction and we conclude. $\square$

### 3.5 Proof of Theorem 4

Now, we have all the tools necessary to prove our main results. We begin with the proof of Theorem 4.

*Proof.* (Proof of Theorem 4) Let $\bar{s}(x,t)$ and $\bar{u}(x,t)$ be the unique entire solutions from Theorem 3. Since $s(x) \equiv 1$ and $u(x) \equiv g(1)$ are supersolutions for (3) and $(0,0)$ are subsolutions, then by standard parabolic theory for monotone systems $s(x,t) \to s(x)$ and $u(x,t) \to u(x)$ uniformly on compact sets, as well as its derivatives, where $s(x)$ and $u(x)$ solve (10) Furthermore, given $\bar{s}_t(x,t), \bar{u}_t(x,t) \geq 0$ then $s(x) \neq 0$ and $u(x) \neq 0$. First, assume that $L \geq L^\star$, from Lemma 4 and Lemma 5 we know that there exist monotonically decreasing steady state solutions

$$\lim_{x \to \infty} s_L(x) = 0 \quad \lim_{x \to \infty} g(s_L(x)) = 0.$$

Furthermore, given

$$\bar{s}(x, -\infty) < s_L(x) \quad \bar{u}(x) < g(s_L(x))$$

then the comparison principle proves $(i)$.

Now, consider the case when $L < L^\star$. In this case, the only steady state solution is the symmetric solution. In fact, from Lemma 8 we know that all steady state solutions must satisfy $s_L(x) \in (b, 1)$ with $s'(x) > 0$ for $x \geq L$. This automatically implies that the steady state solution for the moving average of crime, call it $u_L(x)$ also satisfies $u_L(x) \in (g(b), g(1))$ with $u'_L(x) > 0$ for $x > L$. In fact, we have that

$$\lim_{t \to \infty} s(x,t) = s_m(x) \quad \lim_{t \to \infty} u(x,t) = g(s_m(x)).$$

$\square$

*Remark* 2. Theorem 4 can be extended to solutions $s(x,t)$ and $u(x,t)$ (3) with initial data that vaguely resemble traveling waves, where the fronts of the "waves" are located far enough left of the gap. The solutions will resemble traveling fronts before reaching the gap region.

From proof above we immediately obtain as a corollary the same result for equation (7) with bistable $f(s)$. Recall, that this result was conjectured in [13].

**Corollary 1.** *Let $L^\star$ be as in Theorem 4 and $s(x,t)$ a solution to (7) with $f(s)$ bistable then for the following hold*

(i) *If $L \geq L^\star$ then the wave propagation is blocked.*

(ii) *If $L < L^\star$ then the wave propagates.*

## 3.6 Is Splitting Resources Useful?

In this subsection we focus on the issue of whether splitting resources can improve the situation by minimizing the total resources needed to prevent the propagation. To make this rigorous let us set $L < L^\star$ and arbitrarily split $L$ into two regions of lengths $L_1$ and $L_2$ separated by an interval of length $d$. To summarize we have

$$L_1 + L_2 := L < L^\star.$$

Now, the reaction term, $f(u)$, is replaced by decay term, $-s$, in the two regions of length $L_1$ and $L_2$. The question we aim to answer is this: Is there a distance $d$ between the two gaps of length $L_1$ and $L_2$ so that the wave propagation is prevented? For clarity, we state the problem mathematically by defining a new function.

$$f_d(x,s) = \begin{cases} f(s) & x \in (-\infty, 0) \cap (L_1, L_1+d) \cap (L+d, \infty) \\ -s & x \in [0, L_1] \cap [L_1+d, L+d], \end{cases} \quad (55)$$

and solve the new problem

$$s'' + f_d(x,s) = 0. \quad (56)$$

Figure 2 gives an illustration of the double gap problem. It is self-evident that if the distance between the gaps is very large then, by definition of the critical length, the wave will propagate through the first gap. If the traveling front passes and is allowed to reform, once it hits the second gap the length will once again be too small to stop the propagation of the wave. Hence, we know that for $d$ too large splitting the resources make will necessarily increase the total amount of resources required to prevent the propagation of the waves.

Figure 2: Double Gap Problem

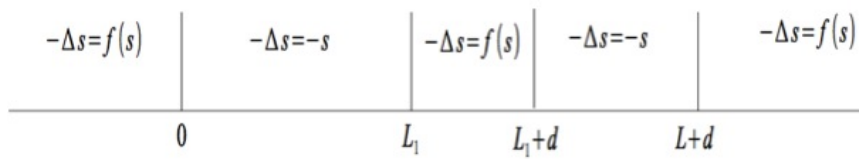

**Lemma 9.** *Let $(s(x), u(x))$ be a solutions to (56) with $L_1 + L_2 = L < L^\star$, then for any $\epsilon > 0$ there exists an $x_\epsilon > L + d$ such that $s(x) > 1 - \epsilon$ and $u(x) > g(1) - \epsilon$ for all $x > x_\epsilon$.*

This lemma implies that splitting resources will never improve the situation, or reduce the total number of resources. Indeed, Lemma (9) proves Proposition 1. In the first step we prove the result for small $d$. The second step is an inductive step, where we prove that if the result holds for an arbitrary $d > 0$, then it holds for $d + \epsilon$ for $\epsilon > 0$.

*Proof.* We prove this in two steps.

*Step 1:* Assume that $0 < d < L^\star - L$ and let $s^\star(x)$ be the strictly decreasing solution obtain in Lemma 5. In particular, $s^\star(x)$ satisfies $s^\star_{xx} + f_{L^\star}(x, s^\star) = 0$. In this case we have $f_{L^\star}(x, s^\star) \leq f_d(x, s^\star)$ which implies that $s^\star$ is a subsolution of

$$s'' + f_d(x, s) = 0.$$

In particular, the inequality is strict as $s^\star(x) \in (b, 1)$ for $x \in [0, L^\star]$. Therefore, $s(x) > s^\star(x)$ and $u(x) \geq g(s^\star(x))$. Thus, $s(L^\star) > b$ which implies, as we saw before, that the wave will propagate and (26) is satisfied.

*Step 2: (Inductive Step)* Assume that for $d = d_1 > 0$ the traveling waves propagate. It suffices to show that separating the gaps by a distance of $d_1 + \epsilon_1$, for $\epsilon_1$ positive, arbitrary, and independent of $d_1$, will not improve the situation. We build a suitable subsolution. Let $s_{d_1}(x)$ be the solution to (56) with $d = d_1$, which by assumption satisfies (26) for some $x_\epsilon > L + d_1$. Furthermore, let $s_s(x)$ be the symmetric solution to

$$\begin{cases} -s'' = -s & L + d_1 \leq x \leq L + d_1 + \epsilon_1 \\ s(L + d_1) = s_{d_1}(L + d_1) & s'(L + d_1) = s'_{d_1}(L + d_1), \end{cases}$$

and $\psi(x)$ be the solution to

$$\begin{cases} -s'' = f(s) & x > L + d_1 + \epsilon_1 \\ s(L + d_1 + \epsilon_1) = s_d(L + d_1). \end{cases}$$

Both $s_s(x)$ and $\psi(x)$ exist and furthermore $\psi(L + d_1 + \epsilon_1) > b$ by definition of $s_{d_1}(x)$. Then

$$\overline{s}(x) = \begin{cases} s_{d_1}(x) & x \leq L + d_1 \\ s_s(x) & L + d_1 \leq x \leq L + d_1 + \epsilon_1 \\ \psi(x) & x > L + d_1 + \epsilon_1, \end{cases}$$

is a subsolution to (56) with $d = d_1 + \epsilon_1$. Indeed, we have that for $x \leq L + d_1$

$$-\overline{s}'' - f_{d_1 + \epsilon_1}(x, \overline{s}) \leq 0,$$

and for $x > L + d_1$ we have that $f_{d_1 + \epsilon_1}(x, s) = f_{d_1 + \epsilon_1}(x, \overline{s})$

$$\overline{s}'' - f_{d_1 + \epsilon_1}(x, \overline{s}) = 0.$$

Hence, $\overline{s}(x)$ is a subsolution and given that $\overline{s}(L + d_1 + \epsilon_1) > b$ then any solution, $s(x)$ must also satisfy $s(L + d_1 + \epsilon_1) > b$. Since the splitting was arbitrary the previous arguments conclude the proof. □

## 4 The Monostable Case: Neutral Criminal Activity

Theorem 5 directly from Theorem 2.2 in [20]. In terms of preventing the invasion of criminal activity, we prove that for all $L > 0$ all steady state solutions increase to one as $x, t \to \infty$. In fact, in this case, due the strength of the instability of zero the solutions will always approach the symmetric solution unless the initial data is exactly zero.

*Proof.* (Proof of Theorem 6) Consider the initial value problem

$$s'' + f(s) = 0 \quad x \in [L, \infty)$$
$$s(L) = \beta,$$

for any $\beta \in (0,1]$. All bounded solutions approach one as $x \to \infty$. As $\lim s(x,t) = \bar{s}(x)$ and $\lim u(x,t) = g(\bar{s}(x))$ on compact sets then it must be the case that $\lim_{x,t \to \infty} s(x,t) = 1$ and $\lim_{t,x \to \infty} u(x,t) = g(1)$.

$\square$

As before we obtain the same result for monostable single parabolic equation.

**Corollary 2.** *Let $L > 0$ be arbitrary and $s(x,t)$ a solution to (7) with $f(s)$ monostable then the wave always propagates.*

## 5 The Single Parabolic Equation

The proof of Theorem 7 follows from Theorem 2.1 in [4] if we study the symmetric problem,

$$s_t - s_{xx} - f_L(x,s) = 0 \tag{57}$$

where,

$$f_L(x,s) = \begin{cases} f(s) & x < -L \text{ and } x > 0 \\ -s & -L \leq x \leq 0. \end{cases} \tag{58}$$

Additionally, one needs to assume that the invasion is from the right. Hence, the traveling wave solution we are interested in satisfies for $z = x + ct$ with $c > 0$

$$\begin{cases} \psi''(z) - c\psi'(z) + f(z) = 0 \\ \psi(-\infty) = 0 \quad \psi(+\infty) = 1. \end{cases} \tag{59}$$

Note that this problem is equivalent (symmetric) to the previous problem. With the help of $\psi(z)$ in the pre-invasion process we show that there is an entire solution on $\mathbb{R} \times \mathbb{R}$, which converge to $\psi(z)$ as $t \to -\infty$. As the proof is similar to that of Theorem 3 and follows the proof of Theorem 2.1 in [4] we leave the details to the reader.

## 6 Discussion

We have studied the effect that a population's *natural tendencies towards crime* have on criminal activity patterns. The simplest case is when there is a natural tendency towards criminal activity, equivalently in our model $s_b > 0$. In this case, in the long run one expects either a constant hotspot or a warm-spot. In the case when there is indifference, $s_b = 0$, lack of criminal activity is a steady state; however, it is unstable and if the payoff is high enough, criminal activity will have a tendency to dominate. The most interesting case is when the population has a natural tendency to be peaceful, or to avoid criminal activity. Here, there is an interesting interplay between the natural tendency and the payoff of committing a crime. In this situation, there are three steady states, complete lack of crime, a small amount of criminal activity, or a hotspot. The lack of criminal activity and hotspots steady state solutions are stable and the small amount of crime is unstable. Hence, in this case it will almost always be either all or nothing, so that in the long run there is a high criminal activity or zero criminal activity.

We focused on the neutral and anti-crime tendencies to study the invasion of criminal activity into areas that start out with zero criminal activity. In both cases, we prove that the invasion of criminal activities is possible, under the right circumstances. This invasion happens via traveling wave solutions. Note that this also holds for $\mathbb{R}^n$ in general, of particular interest in this application is the $\mathbb{R}^2$ case. Indeed in Figure 3 we see a traveling wave solution for $s(x,t)$ moving in the $x$-direction. In the neutral case, or monostable case as we have been referring to it, the high criminal activity steady state solution is always more stable and the invasion will always be that direction. On the other hand, in the anti-criminal activity case, or bistable case, it is easier to have the high crime steady state solution be more stable; however, there is the possibility that the invasion to be in the other direction. More specifically this implies that there are parameters that lead to zero crime areas invading the high crime rate zones.

Figure 3: 2D Traveling Wave of the Propensity to Commit a Crime

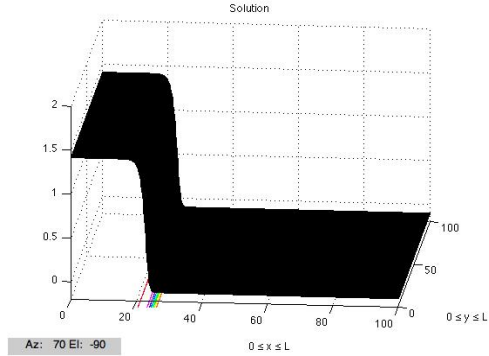

We have also addressed the issue of the prevention of the propagation of criminal activity by applying a finite amount of resources to reduce the payoff of committing a crime to zero on a finite interval. This corresponds to employing crime deterring strategies such a providing more police officers to the area, decreasing opportunities, etc. In the bistable case, we showed that there is a critical length of the interval with zero payoff coefficient, or equivalently a minimum amount of resources, required to prevent the propagation of criminal activity. The amount of resources required depends on $\mathcal{F}(1)$, vaguely measures the stability of the hotspot steady state, as was shown numerically. An interesting conclusion of our work is that in the monostable case there is no possible amount of *finite* resources that will stop the propagation of crime. Hence, to stop the propagation of crime it is essential that the population be naturally anti-criminal activity or to start out with zero-crime.

There are many interesting questions left to answer. For example, in the anti-criminal activity scenario there are two ways to prevent the invasion of the criminal activity waves. The first is to have coefficients such that $\int f(s)\,ds < 0$ and the second is to employ enough resources on an area of length $L^\star$ to diminish the payoff to zero there. An interesting question is to determine which of these options is more cost effective. In the former case corresponds to changing both the payoff to commit a crime and the natural crime-tendencies, this could be done say via preventive educational efforts. The latter is simply a direct intervention in reducing the payoff of committing a crime.

## A  Lemma 2

We have shown that there exists traveling wave solutions, $(\psi(z), \phi(z))$, with speed $c > 0$ which satisfy

$$\psi''(z) + c\psi'(z) - \psi + \alpha\phi(z) = 0 \tag{60a}$$
$$c\phi'(z) + g(\psi(z)) - \phi(z) = 0 \tag{60b}$$

This can be turned into a system of ODE's by defining $(v, q, w) = (\psi, \phi, \psi')$. Then the above system is equivalent to

$$v' = w$$
$$q' = \frac{q - g(v)}{c}$$
$$w' = -cw + v - \alpha q.$$

The Jacobian at a steady state $\mathbf{v}^* = (v^*, q^*, w^*)$ is

$$J(\mathbf{v}^*) = \begin{bmatrix} 0 & 0 & 1 \\ \frac{-g'(v^*)}{c} & \frac{1}{c} & 0 \\ 1 & -\alpha & -c \end{bmatrix},$$

with the characteristic polynomial

$$\left(\lambda - \frac{1}{c}\right)\left(\lambda^2 + c\lambda - 1\right) - \frac{\alpha g'(v^*)}{c}.$$

Recall that $g'(0) = 0$ and so for the steady state $(0, 0)$ the roots of the characteristic polynomial are

$$\lambda_1 = \frac{1}{c}, \ \lambda_\pm = \frac{-c \pm \sqrt{c^2 + 4}}{2}.$$

In fact, $\lambda = \lambda_-$. Note that the characteristic polynomial is increasing at $\lambda = \lambda_-$. For $k \in \mathbb{R}$ we now define the functional

$$Q(\lambda, k) := \left(\lambda - \frac{1}{c}\right)\left(\lambda^2 + c\lambda - 1\right) - k. \tag{61}$$

As $k = -\frac{g(1)\alpha}{c} < 0$ we see that the negative root of the polynomial must decrease and from this we conclude that $\lambda > \mu$.